\documentclass{article}
\usepackage{amsmath,amssymb}

\author{Gert Almkvist}
\title{Calabi-Yau differential equations of degree 2 and 3 and Yifan
  Yang's pullback}

\begin{document}

\maketitle
\textbf{1.Introduction.}

The start of [2] was the pullback of a 5-th order differential equation,
coming from diophantine approximation in number theory, to a 4-th order
differential equation of Calabi-Yau type. Recently Yifan Yang has found a
different way to do this pullback. His motivation came from a wish to
simplify the matrices generating the monodromy group. Here we get the same
pullback in a more direct way.

There are 14 hypergeometric Calabi-Yau ( called CY in the sequel)
differential equations (see [1],[3]), i.e. of type 
\[
\theta ^{4}+xP_{1}(\theta )
\]%
where $\theta =x\frac{d}{dx}$ and degree $P_{1}=4$. The next step is CY
equations of degree $2$, i.e. of type 
\[
\theta ^{4}+xP_{1}(\theta )+x^{2}P_{2}(\theta )
\]%
where degree $P_{1},P_{2}\leq 4.$ So far we know 67 CY equations of degree $2
$. There are $4\cdot 10-2=38$ Hadamard products of two second order
equations of very special type (see [2] and section 3.1). We find 8 more
coming from third order equations by multiplying the coefficients by $\binom{%
2n}{n}.$Then there are $14$ Yifan Yang pullbacks of special $5$-th order
hypergeometric equations. Finally there are $7$ sporadic equations found
through a computer search by van Enckevort and van Straten (see [5]).

One can continue to ask for CY equations of degree $3.$ So far we know only $%
7$ such examples. Observe that we count only equations of minimal degree
under the equivalence relation having the same instanton numbers.

In section 4 we consider difference equations of order $2$ and $3$ (usually
of degree $6$ ), leading to CY equations after factorization in the Weyl
algebra.

In the course of defining the Yifan Yang pullback we find a number of
formulas for various Wronskians of the $4$ solutions of a CY equation. 
\[
\]

\textbf{2. Wronskians and Yifan Yang's pullback.}

\[
\]
\textbf{2.1 The double Wronskian is almost the square.}

Consider the 4-th order differential equation 
\[
y^{^{\prime \prime \prime \prime \prime }}+a_{3}y^{\prime \prime \prime
}+a_{2}y^{\prime \prime }+a_{1}y^{\prime }+a_{0}y=0
\]%
\[
\]%
We assume it is MUM (Maximal Unipotent Monodromy) so we have the Frobenius
solutions 
\[
y_{0}=1+A_{1}x+...
\]%
\[
y_{1}=y_{0}\log (x)+B_{1}x+..
\]%
\[
y_{2}=\frac{1}{2}y_{0}\log (x)^{2}+...
\]%
\[
y_{3}=\frac{1}{6}y_{0}\log (x)^{3}+...
\]%
\[
\]%
We also assume that the differential equation is Calabi-Yau (C-Y2), i.e. the
coefficients satisfy 
\[
a_{1}=\frac{1}{2}a_{2}a_{3}-\frac{1}{8}a_{3}^{3}+a_{2}^{\prime }-\frac{3}{4}%
a_{3}a_{3}^{\prime }-\frac{1}{2}a_{3}^{\prime \prime }
\]%
\[
\]

Let 
\[
t=\frac{y_1}{y_0} 
\]
Then in [1] it is proved that the condition C-Y2 is equivalent to 
\[
\frac{d^2}{dt^2}(\frac{y_2}{y_0})=\frac 1{y_0^2(\frac{dt}{dx})^3}\exp
(-\frac 12\int a_3dx) 
\]
\[
\]
In [2] it is proved that C-Y2 impies that 
\[
w_0=x\left| 
\begin{array}{rr}
y_0 & y_1 \\ 
y_0^{\prime } & y_1^{\prime }%
\end{array}
\right| 
\]
\[
w_1=x\left| 
\begin{array}{rr}
y_0 & y_2 \\ 
y_0^{\prime } & y_2^{\prime }%
\end{array}
\right| 
\]
\[
w_2=x\left| 
\begin{array}{rr}
y_0 & y_3 \\ 
y_0^{\prime } & y_3^{\prime }%
\end{array}
\right| =x\left| 
\begin{array}{rr}
y_1 & y_2 \\ 
y_1^{\prime } & y_2^{\prime }%
\end{array}
\right| 
\]
\[
w_3=\frac x2\left| 
\begin{array}{rr}
y_1 & y_3 \\ 
y_1^{\prime } & y_3^{\prime }%
\end{array}
\right| 
\]
\[
w_4=\frac x2\left| 
\begin{array}{rr}
y_2 & y_3 \\ 
y_2^{\prime } & y_3^{\prime }%
\end{array}
\right| 
\]
\[
\]
form a Frobenius basis of the solutions to a 5-th order differential
equation 
\[
w^{(5)}+b_4w^{(4)}+b_3w^{\prime \prime \prime }+b_2w^{\prime \prime
}+b_1w^{\prime }+b_0w=0 
\]
\[
\]
where 
\[
a_3=\frac 25b_4+\frac 2x 
\]
Now consider the double Wronskian 
\[
u=\left| 
\begin{array}{rr}
w_0 & w_1 \\ 
w_0^{\prime } & w_1^{\prime }%
\end{array}
\right| 
\]
\[
\]
Differentiating $t=\dfrac{y_1}{y_0}$ we get 
\[
\frac{dt}{dx}=\frac{y_0y_1^{\prime }-y_0^{\prime }y_1}{y_0^2}=\frac{w_0}{%
xy_0^2} 
\]
Similarly 
\[
\frac d{dt}(\frac{y_2}{y_0})=(\frac{y_2}{y_0})^{\prime }\cdot \frac{dx}{dt}=%
\frac{(y_0y_2^{\prime }-y_0^{\prime }y_2)xy_0^2}{y_0^2w_0}=\frac{w_1}{w_0} 
\]
Then 
\[
\frac{d^2}{dt^2}(\frac{y_2}{y_0})=(\frac{w_1}{w_0})^{\prime }\cdot \frac{dx}{%
dt}=\frac{xy_0^2(w_0w_1^{\prime }-w_0^{\prime }w_1)}{w_0^3}=\frac{xy_0^2u}{%
w_0^3} 
\]
But by the condition equivalent to C-Y2 this is equal to 
\[
\frac 1{y_0^2(\dfrac{dt}{dx})^3}\exp (-\frac 12\int a_3dx) 
\]
Hence 
\[
u=\left| 
\begin{array}{rr}
w_0 & w_1 \\ 
w_0^{\prime } & w_1^{\prime }%
\end{array}
\right| =x^2y_0^2\exp (-\frac 12\int a_3dx) 
\]
\[
\]

Yifan Yang has found the following pullback of the 5-th order equation 
\[
Y_0=x^{5/2}u^{1/2}\exp (-\frac 15\int b_4dx) 
\]
which is 
\[
Y_0=x^{5/2}xy_0\exp (-\int (\frac 14a_3+\frac 15b_4)dx) 
\]
\[
=x^{7/2}y_0\exp (-\int (\frac 14(\frac 25b_4+\frac 2x)+\frac 15b_4)dx) 
\]
\[
=x^3y_0\exp (-\frac 3{10}\int b_4dx) 
\]
\[
\]

Hence Yifan Yang's pullback is just $y_0$ multiplied by a suitable factor.
The remarkable fact is that for most of our 5-th order differential
equations his pullback satisfies a 4-th order equation where the degree is
half of the degree of the equation for $y_0.$ This is not the case for a
general $y_0.$ On the contrary, the substitution 
\[
Y_0=x^{3/2}y_0\exp (-\frac 34\int a_3dx) 
\]
\[
\]
of e.g. \#209 gives an equation of degree 16 (\#209 refers to the numbering
in [3] also called the big table)

\textbf{Example. }(\#130 by Helena Verrill)

Let 
\[
A_n^{\prime }=\sum_{i+j+k+l+m+p=n}\left( \frac{n!}{i!j!k!l!m!p!}\right) ^2 
\]
Then 
\[
w_0=\sum A_n^{\prime }x^n 
\]
\[
\]
satisfies the 5-th order differential equation 
\[
\theta ^5-2x(2\theta +1)(14\theta ^4+28\theta ^3+28\theta ^2+14\theta +3) 
\]
\[
+4x^2(\theta +1)^3(196\theta ^2+392\theta +255)-1152x^3(\theta +1)^2(\theta
+2)^2(2\theta +3) 
\]
\[
\]
The Yifan pullback is 
\[
\theta ^4-4x\left\{ 28(\theta +\frac 12)^4+28(\theta +\frac 12)^2+1\right\} 
\]
\[
+3x^2\left\{ 1568(\theta +1)^4+2130(\theta +1)^2+225\right\} 
\]
\[
-4x^3\left\{ 23104(\theta +\frac 32)^432532(\theta +\frac 32)^2+3213\right\} 
\]
\[
+x^4\left\{ 872704(\theta +2)^4+995680(\theta +2)^2+93337\right\} 
\]
\[
-2^93^2x^5(\theta +\frac 52)^2\left\{ 784(\theta +\frac 52)^2+647\right\} 
\]
\[
+2^{14}3^4x^6(2\theta +1)^2(\theta +\frac 52)(\theta +\frac 72) 
\]
which is of degree 6 while the ordinary pullback in the big table is of
degree 12) 
\[
\]

\textbf{2.2 Relations between the bases of the Wronskians of the 5-th order
equation and the symmetric squares of the 4-th order equation satisfying
C-Y2.}

The space of the Wronskians of the solutions of the 5-th order equation has
dimension $\binom 52=10$. Similarly is the dimension of the space of
symmetric squares of the solutions of the 4-th order equation also $\binom{%
4+1}2=10$ (pointed out by van Straten).

Let 
\[
f:=\exp (-\frac 12\int a_3dx) 
\]
Then we have the following explicit formulas 
\[
\left| 
\begin{array}{rr}
w_0 & w_1 \\ 
w_0^{\prime } & w_1^{\prime }%
\end{array}
\right| =x^2f\cdot y_0^2 
\]
\[
\left| 
\begin{array}{rr}
w_0 & w_2 \\ 
w_0^{\prime } & w_2^{\prime }%
\end{array}
\right| =x^2f\cdot y_0y_1 
\]
\[
\left| 
\begin{array}{rr}
w_0 & w_3 \\ 
w_0^{\prime } & w_3^{\prime }%
\end{array}
\right| =\frac 12x^2f\cdot y_1^2 
\]
\[
\left| 
\begin{array}{rr}
w_0 & w_4 \\ 
w_0^{\prime } & w_4^{\prime }%
\end{array}
\right| =\frac 12x^2f\cdot (y_1y_2-y_0y_3) 
\]
\[
\left| 
\begin{array}{rr}
w_1 & w_2 \\ 
w_1^{\prime } & w_2^{\prime }%
\end{array}
\right| =x^2f\cdot y_0y_2 
\]
\[
\left| 
\begin{array}{rr}
w_1 & w_3 \\ 
w_1^{\prime } & w_3^{\prime }%
\end{array}
\right| =\frac 12x^2f\cdot (y_1y_2+y_0y_3) 
\]
\[
\left| 
\begin{array}{rr}
w_1 & w_4 \\ 
w_1^{\prime } & w_4^{\prime }%
\end{array}
\right| =\frac 12x^2f\cdot y_2^2 
\]
\[
\left| 
\begin{array}{rr}
w_2 & w_3 \\ 
w_2^{\prime } & w_3^{\prime }%
\end{array}
\right| =\frac 12x^2f\cdot y_1y_3 
\]
\[
\left| 
\begin{array}{rr}
w_2 & w_4 \\ 
w_2^{\prime } & w_4^{\prime }%
\end{array}
\right| =\frac 12x^2f\cdot y_2y_3 
\]
\[
\left| 
\begin{array}{rr}
w_3 & w_4 \\ 
w_3^{\prime } & w_4^{\prime }%
\end{array}
\right| =\frac 14x^2f\cdot y_3^2 
\]
\[
\]

\textbf{2.3 Higher order Wronskians.}

We consider also Wronskians of order 3 and 4. Let in general 
\[
W(u_1.u_2,u_3)=\left| 
\begin{array}{rrr}
u_1 & u_2 & u_3 \\ 
u_1^{\prime } & u_2^{\prime } & u_3^{\prime } \\ 
u_1^{\prime \prime } & u_2^{\prime \prime } & u_3^{\prime \prime }%
\end{array}
\right| 
\]
\[
\]
and similarly for 4-th order Wronskians. 
\[
\]

\textbf{Proposition: }We have the following identities 
\[
W((y_0,y_1,y_2)=fy_0 
\]
\[
W(y_0,y_1,y_3)=fy_1 
\]
\[
W(y_0,y_2,y_3)=fy_2 
\]
\[
W(y_1,y_2,y_3)=fy_3 
\]

\textbf{Proof.} We have earlier proved 
\[
W(w_0,w_1)=x^2fy_0^2 
\]
But a direct calculation shows that 
\[
W(w_0,w_1)=x^2y_0W(y_0,y_1,y_2) 
\]
The other identities are proved similarly. 
\[
\]

\textbf{Corollary: }$W(y_0,y_1,y_2),$ etc satisfy a 4-th order differential
equation with the same mirror map and Yukawa coupling as the original
equation. Usually the degree of the equation is doubled. 
\[
\]

We have the following identities 
\[
W(w_0,w_1,w_2)=x^3W(y_0.y_1,y_2)^2=x^3f^2y_0^2 
\]
\[
W(w_0,w_1,w_3)=x^3W(y_0,y_1,y_2)W(y_0,y_1,y_3)=x^3f^2y_0y_1 
\]
\[
W(w_0,w_1,w_4)=x^3W(y_0,y_1,y_2)W(y_0,y_2,y_3)=x^3f^2y_0y_2 
\]
\[
W(w_0,w_2,w_3)=\frac 12x^3W(y_0,y_1,y_3)^2=x^3f^2y_1^2 
\]
\[
W(w_0,w_2,w_4)=\frac 12x^3\left\{
W(y_0,y_1.y_2)W(y_1,y_2,y_3)+W(y_0,y_1,y_3)W(y_0,y_2,y_3)\right\} 
\]
\[
=\frac 12x^3f^2(y_0y_3+y_1y_2) 
\]
\[
W(w_0,w_3,w_4)=\frac 12x^3W(y_0,y_1,y_3)W(y_1,y_2,y_3)=\frac 12x^3f^2y_1y_3 
\]
\[
W(w_1,w_2,w_3)=\frac 12x^3\left\{
W(y_0,y_1,y_3)W(y_0,y_2,y_3)-W(y_0,y_1,y_2)W(y_1,y_2,y_3)\right\} 
\]
\[
=\frac 12x^3f^2(y_1y_2-y_0y_3) 
\]
\[
W(w_1,w_2,w_4)=\frac 12x^3W(y_0,y_2,y_3)^2=\frac 12x^3f^2y_2^2 
\]
\[
W(w_1,w_3,w_4)=\frac 12x^3W(y_0,y_2,y_3)W(y_1,y_2,y_3)=\frac 12x^3f^2y_2y_3 
\]
\[
W(w_2,w_3,w_4)=\frac 14x^3W(y_1,y_2,y_3)^2=\frac 14x^3f^2y_3^2 
\]
\[
\]
We have the following formulas for 4-th order Wronskians 
\[
W(w_0,w_1,w_2,w_3)=x^3f^3w_0 
\]
\[
W(w_0,w_1,w_2,w_4)=x^3f^3w_1 
\]
\[
W(w_0,w_1,w_3,w_4)=x^3f^3w_2 
\]
\[
W(w_0,w_2,w_3,w_4)=x^3f^3w_3 
\]
\[
W(w_1,w_2,w_3,w_4)==x^3f^3w_4 
\]
\[
\]
We emphasize that all these formulas are only valid under the condition C-Y2.

\[
\]
\textbf{2.4 The secret of the Yifan Yang pullback.}

\smallskip Consider the differential equation 
\[
y^{(4)}+a_3y^{\prime \prime \prime }+a_2y^{\prime \prime }+a_1y^{\prime
}+a_0y=0 
\]
which we assume is MUM and satisfying the C-Y2 condition 
\[
a_1=\frac 12a_2a_3-\frac 18a_3^3+a_2^{\prime }-\frac 34a_3a_3^{\prime
}-\frac 12a_3^{\prime \prime } 
\]
Then the wronskian 
\[
w_0=x\left| 
\begin{array}{rr}
y_0 & y_1 \\ 
y_0^{\prime } & y_1^{\prime }%
\end{array}
\right| 
\]
satisfies the 5-th order equation 
\[
w^{(5)}+b_4w^{(4)}+b_3w^{\prime \prime \prime }+b_2w^{\prime \prime
}+b_1w^{\prime }+b_0w=0 
\]
We have the following relations between the coefficients 
\[
a_3=2x^{-1}+\frac 25b_4 
\]
\[
a_2=-\frac 32x^{-2}+\frac 12b_3+\frac 35x^{-1}b_4-\frac 7{50}b_4^2-\frac
25b_4^{\prime } 
\]
\[
a_1=\frac 32x^{-3}+b_2+\frac 12x^{-1}b_3-b_3^{\prime }-\frac 12b_3b_4-\frac
3{10}x^{-2}b_4 
\]
\[
-\frac 7{50}x^{-1}b_4^2+\frac{31}{250}b_4^3-\frac 25x^{-1}b_4^{\prime }+%
\frac{18}{25}b_4b_4^{\prime } 
\]
\[
a_0=-\frac{15}{16}x^{-4}-\frac 14b_1+\frac 14x^{-1}b_2+\frac 12b_2^{\prime
}+\frac 3{20}b_2b_4 
\]
\[
-\frac 18x^{-2}b_3+\frac 1{16}b_3^2-\frac 18x^{-1}b_3^{\prime }-\frac
38b_3^{\prime \prime }-\frac 1{10}x^{-1}b_3b_4 
\]
\[
-\frac{21}{200}b_3b_4^2-\frac 7{20}b_3b_3^{\prime }-\frac{13}{40}b_3^{\prime
}b_4+\frac 3{20}x^{-3}b_4+\frac 7{200}x^{-2}b_4^2 
\]
\[
+\frac{11}{500}x^{-1}b_4^3+\frac{221}{10000}b_4^4+\frac
1{10}x^{-2}b_4^{\prime }-\frac 1{20}x^{-1}b_4^{\prime \prime }+\frac
1{10}b_4^{\prime \prime \prime } 
\]
\[
+\frac 3{50}x^{-1}b_4b_4^{\prime }+\frac{29}{100}b_4b_4^{\prime \prime }+%
\frac{33}{100}(b_4^{\prime })^2+\frac{69}{250}b_4^2b_4^{\prime } 
\]
\[
\]
We express the condition C-Y2 in the coefficients $b_2,b_3,b_4$%
\[
U\equiv -b_2+\frac 32b_3^{\prime }+\frac 35b_3b_4-b_4^{\prime \prime }-\frac
65b_4b_4^{\prime }-\frac 4{25}b_4^3=0 
\]
\[
\]
Now we assume the 5-th order equation given satisfying the condition just
stated and we want to find the 4-th order equation somewhat modified to get
lower degrees of the coefficients. To achieve this we will make the
transformation 
\[
u=f\cdot y 
\]
Then the 4-th order equation is transformed to 
\[
u^{(4)}+c_3u^{\prime \prime \prime }+c_2u^{\prime \prime }+c_1u^{\prime
}+c_0u=0 
\]
where 
\[
c_3=4\frac{f^{\prime }}f+a_3 
\]
\[
c_2=6\frac{f^{\prime \prime }}f+3a_3\frac{f^{\prime }}f+a_2 
\]
\[
c_1=4\frac{f^{\prime \prime \prime }}f+3a_3\frac{f^{\prime \prime }}f+2a_2%
\frac{f^{\prime }}f+a_1 
\]
\[
c_0=\frac{f^{(4)}}f+a_3\frac{f^{\prime \prime \prime }}f+a_2\frac{f^{\prime
\prime }}f+a_1\frac{f^{\prime }}f+a_0 
\]
\[
\]
To simplify 
\[
a_3=\frac 2x+\frac 25b_4 
\]
we choose 
\[
f=x^{-1/2}\exp (\lambda \int b_4dx) 
\]
with 
\[
\frac{f^{\prime }}f=-\frac 1{2x}+\lambda b_4 
\]
\[
\frac{f^{\prime \prime }}f=\frac 34x^{-2}-\lambda x^{-1}b_4+\lambda
^2b_4^2+\lambda b_4^{\prime } 
\]
\[
\frac{f^{\prime \prime \prime }}f=-\frac{15}8x^{-3}+\frac 94\lambda
x^{-2}b_4-\frac 32\lambda ^2x^{-1}b_4^2+\lambda ^3b_4^3-\frac 32\lambda
x^{-1}b_4^{\prime }+\lambda b_4^{\prime \prime } 
\]
\[
\frac{f^{(4)}}f=\frac{105}{16}x^{-4}-\frac{15}2\lambda x^{-3}b_4+\frac
92\lambda ^2x^{-2}b_4^2-2\lambda ^3x^{-1}b_4^3+\lambda ^4b_4^4 
\]
\[
+\frac 92\lambda x^{-2}b_4^{\prime }-2\lambda x^{-1}b_4^{\prime \prime
}+\lambda b_4^{\prime \prime \prime }+4\lambda ^2b_4b_4^{\prime \prime
}+6\lambda ^3b_4^2b_4^{\prime }+3\lambda ^2(b_4^{\prime })^2 
\]
\[
\]
We get rid of the term $\frac 2x$ in $a_3$%
\[
c_3=\frac{2(10\lambda +1)}5b_4 
\]
\[
c_2=\frac 12b_3+(6\lambda ^2+\frac 65\lambda -\frac 7{50})b_4^2+(6\lambda
-\frac 25)b_4^{\prime } 
\]
\[
\]
We want to to avoid higher powers and derivatives of $b_4$ (remember that
the $b_k$ are rational functions so differentiating also increases the
degree). We can also add multiples of the relation $U$ (and its
derivatives). Hence consider 
\[
c_1+\mu U= 
\]
\[
(1-\mu )b_2+(-1+\frac 32\mu )b_3^{\prime }+(\lambda +\frac 35\mu -\frac
12)b_3b_4 
\]
\[
+(4\lambda ^3+\frac 65\lambda ^2-\frac 75\lambda +\frac{31}{250}-\frac
4{25}\mu )b_4^3 
\]
\[
+(12\lambda ^2+\frac 25\lambda +\frac{18}{25}-\frac 65\mu )b_4b_4^{\prime } 
\]
\[
+(4\lambda -\mu +\frac 25)b_4^{\prime \prime } 
\]
\[
\]
By choosing $\lambda $ and $\mu $ such that the three last terms vanish we
get a very simple formula for $c_1$ agreeing with Yifan's formula 
\[
4\lambda ^3+\frac 65\lambda ^2-\frac 75\lambda +\frac{31}{250}-\frac
4{25}\mu =0 
\]
\[
12\lambda ^2+\frac 25\lambda +\frac{18}{25}-\frac 65\mu =0 
\]
\[
4\lambda -\mu +\frac 25=0 
\]
\[
\]
with the solution 
\[
\lambda =\frac 3{10}\text{ } 
\]
\[
\mu =\frac 85 
\]
giving 
\[
c_1=-\frac 35b_2+\frac 75b_3^{\prime }+\frac{19}{25}b_3b_4 
\]
Fortunately $\lambda =\frac 3{10}$ works also for $c_0$ . Indeed 
\[
c_0+\frac 25U^{\prime }-\frac 1{4x}U+\frac{41}{100}b_4U= 
\]
\[
-\frac 14b_1+\frac 1{10}b_2^{\prime }+\frac 1{25}b_2b_4+\frac
9{40}b_3^{\prime \prime }+\frac 1{16}b_3^2+\frac 1{25}b_3b_4^{\prime } 
\]
\[
+\frac{23}{100}b_3^{\prime }b_4+\frac 9{250}b_3b_4^2 
\]
\[
\]
which agrees with Yifan Yang's formula. To get the analytic solution we have
to multiply the solution by $x^{5/2}$ which is the same thing as the
substitution $\theta \longrightarrow \theta -\frac 52$ in the differential
equation. 
\[
\]

\textbf{Definition,} The \textsl{Yifan Yang pullback }of the 5-th order
equation 
\[
y^{(5)}+b_4y^{(4)}+b_3y^{\prime \prime \prime }+b_2y^{\prime \prime
}+b_1y^{\prime }+b_0y=0 
\]
is 
\[
y^{(4)}+c_3y^{\prime \prime \prime }+c_2y^{\prime \prime }+c_1y^{\prime
}+c_0y=0 
\]
where 
\[
c_3=\frac 85b_4 
\]
\[
c_2=\frac 12b_3+\frac 75b_4^{\prime }+\frac{19}{25}b_4^2 
\]
\[
c_1=-\frac 35b_2+\frac 75b_3^{\prime }+\frac{19}{25}b_3b_4 
\]
\[
c_0=-\frac 14b_1+\frac 1{10}b_2^{\prime }+\frac 1{25}b_2b_4+\frac
9{40}b_3^{\prime \prime }+\frac 1{16}b_3^2 
\]
\[
+\frac 1{25}b_3b_4+\frac{23}{100}b_3^{\prime }b_4+\frac 9{250}b_3b_4^2 
\]
\[
\]
\[
\]

\textbf{3.Calabi-Yau differential equations of degree 2 and 3.}

In the big table [3] of CY equations there are now 336 entries. Except for
the first 16 there is no order in the table. The entries were put in as they
were found. Here we will start to bring some order into the chaos.

Given a (MUM) differential operator L in standard form 
\[
L=\theta ^k+xP_1(\theta )+...+x^dP_d(\theta ) 
\]
where deg($P_j)\leq k$ we call $k$ the \emph{order} and $d$
the \textsl{degree }of $L.$ For a CY equation of degree $d$%
\[
L=\theta ^4+xP_1(\theta )+..+x^dP_d(\theta ) 
\]
having Frobenius solutions $y_0,y_1,y_2,y_3$ we define 
\[
q=\exp (y_1/y_0) 
\]
\[
K(q)=(q\frac d{dq})^2(y_2/y_0)=1+\sum_{k=1}^\infty \frac{k^3N_kq^k}{1-q^k} 
\]
the \textsl{Yukawa coupling.} The $N_k$ are called \textsl{instanton
numbers. }One of the conditions for being a CY equation is that there should
exist a fix integer $c$ such that $cN_k$ are integers for all $k.$. It is
clear that $q$ and $K$ are invariant under the transformation 
\[
y(x)\rightarrow f(x)y(x) 
\]
One can verify that $K$ (but not $q)$ is also invariant under 
\[
y(x)\rightarrow y(g(x)) 
\]
where 
\[
g(x)=x+c_2x^2+c_3x^3+.... 
\]
We call two CY differential equations \textsl{equivalent} ($L_1\sim L_2)$ if
the have the same $K(q)$. So $K(q)$ is the invariant we shall use to
classify CY equations. This is also done in the ''Superseeker'' in the big
table, where the equations are ordered after the absolute values if $N_1$
and $N_3$ (observe that the instanton numbers are not invariant under the
transformation $x\rightarrow cx$ if $c\neq 1$ ). There are many equivalent
equations in the table. Now we will only give the representative of lowest
degree. Thus \#129 in the table is of degree $24$ but it is equivalent to \#$%
\widehat{8}$ , which in its turn will be shown to be equivalent to $%
\widetilde{8}$, a Yifan Yang pullback of degree $2.$

There are equivalences between various Hadamard products, for some of which
we can even find explicit transformations. This will be treated at the end
of this section.

\[
\]

\textbf{3.1 Hadamard products.}

We start with second order equations. They are of two types

\textbf{(i) Hypergeometric.}

They are of type 
\[
\theta ^2-xQ(\theta ) 
\]
where $Q(\theta )$ is a polynomial of degree $2$%
\[
\begin{tabular}{|c|c|}
\hline
Name & $Q(\theta )$ \\ \hline
A & $4(2\theta +1)^2$ \\ \hline
B & $3(3\theta +1)(3\theta +2)$ \\ \hline
C & $4(4\theta +1)(4\theta +3)$ \\ \hline
D & $12(6\theta +1)(6\theta +5)$ \\ \hline
\end{tabular}
\]
\[
\]

\textbf{(ii) Second degree}

They are of type 
\[
\theta ^2-xP(\theta )-cx^2(\theta +1)^2 
\]
where 
\[
P(\theta )=a\theta ^2+a\theta +b 
\]
\[
\begin{tabular}{|c|c|c|c|}
\hline
Name & a & b & c \\ \hline
(a) & 7 & 2 & 8 \\ \hline
(b) & 11 & 3 & 1 \\ \hline
(c) & 10 & 3 & -9 \\ \hline
(d) & 12 & 4 & -32 \\ \hline
(e) & 32 & 12 & -256 \\ \hline
(f) & 9 & 3 & -27 \\ \hline
(g) & 17 & 6 & -72 \\ \hline
(h) & 54 & 21 & -729 \\ \hline
(i) & 128 & 52 & -4096 \\ \hline
(j) & 864 & 372 & -186624 \\ \hline
\end{tabular}
\]
\[
\]

Most of these differential equations appear in Zagier's secret paper [7] and
the rest were found by van Enckevort and Van Straten in a computer search.
More information can be found in [2] where formulas for the coefficients are
given.The Hadamard product is defined by 
\[
\sum A_nx^n*\sum B_nx^n=\sum A_nB_nx^n 
\]
We use $L_1*L_2$ for the corresponding differential operator annihilating $%
\sum A_nB_nx^n$ Then we have the formula 
\[
(\theta ^2-xP(\theta )-cx^2(\theta +1)^2)*(\theta ^2-xQ(\theta )) 
\]
\[
=\theta ^4-xP(\theta )Q(\theta )-cx^2Q(\theta )Q(\theta +1) 
\]
\[
\]
\[
\begin{tabular}{|c|c|c|c|c|}
\hline
{*} & A & B & C & D \\ \hline
a & 45 & 15 & 68 & 62 \\ \hline
b & 25 & 24 & 51 & 63 \\ \hline
c & 58 & 70 & 69 & 64 \\ \hline
d & 36 & 48 & 38 & 65 \\ \hline
e & 111=10$^{**}\thicksim $10$^{*}$ & 110 & 30$\thicksim $3 & 112 \\ \hline
f & 133 & 134 & 135 & 136 \\ \hline
g & 137 & 138 & 139 & 140 \\ \hline
h & 141=8$^{**}$=8$^{*}$ & 142 & K(q)=1 & 143 \\ \hline
i & 7$^{**}\thicksim $7$^{*}$ & B*(i) & C*(i)$\thicksim $6$^{*}$ & D*(i) \\ 
\hline
j & 9$^{**}\thicksim $9$^{*}$ & B*(j) & C*(j) & D*(j) \\ \hline
\end{tabular}
\]
\[
\]
We get $40$ degree $2$ CY equations $(A)*(a),(A)*(b),...,(D)*(j)$ but two of
them,are deleted, since $(C)*(h)$ has trivial $K(q).$and $C*e\thicksim 3.$
We compute the last equation.

\textbf{Example.}$(D)*(j)$

We have for (D) 
\[
Q(\theta )=12(6\theta +1)(6\theta +5) 
\]
where 
\[
\theta ^2-xQ(\theta ) 
\]
has solution 
\[
y_0=\sum \frac{(6n)!}{(3n)!(2n)!n!}x^n 
\]
We have for (j) 
\[
\theta ^2-12x(72\theta ^2+72\theta +31)+186624x^2(\theta +1)^2 
\]
with solution 
\[
y_0=\sum_{n=0}^\infty 432^n\sum_{k=0}^n(-1)^k\binom{-5/6}k\binom{-1/6}{n-k}%
^2x^n 
\]
We get the Hadamard product (D)*(j) 
\[
\theta ^4-144x(6\theta +1)(6\theta +5)(72\theta ^2+72\theta +31) 
\]
\[
+12^2186624x^2(6\theta +1)(6\theta +5)(6\theta +7)(6\theta +11) 
\]
with the solution 
\[
y_0=\sum A_nx^n 
\]
where 
\[
A_n=\frac{(6n)!}{((3n)!(2n)!n1}\sum_{k=0}^n(-1)^k\binom{-5/6}k\binom{-1/6}{%
n-k}^2 
\]
\[
\]
There are also 10 third order equations $\alpha ,\beta ,...,\varkappa .$
(see [2]). Multiplying the coefficients by $\binom{2n}n$ we get the
following 4-th order equations 
\[
\theta ^4-x(2\theta +1)^2P(\theta )+cx^2(\theta +1)^2(2\theta +1)(2\theta
+3) 
\]
where the polynomial $P(n)$ and $c$ are given in the following table 
\[
\]
\[
\begin{tabular}{|c|c|c|c|}
\hline
Name & \# & P(n) & c \\ \hline
$\alpha $ & $16$ & $4(5n^2+5n+2)$ & $256$ \\ \hline
$\gamma $ & $29$ & $16(2n^2+2n+1)$ & $1024$ \\ \hline
$\delta $ & $41$ & $2(7n^2+7n+3)$ & $324$ \\ \hline
$\varepsilon $ & $42$ & $8(3n^2+3n+1)$ & $64$ \\ \hline
$\zeta $ & $185$ & $6(3n^2+3n+1)$ & $-108$ \\ \hline
$\eta $ & $184$ & $2(11n^2+11n+5)$ & $500$ \\ \hline
$\iota $ & $4^{*}$ & $6(9n^2+9n+5)$ & $2916$ \\ \hline
$\varkappa $ & $13^{*}$ & $48(18n^2+18n+13$ & $746496$ \\ \hline
\end{tabular}
\]
\[
\]
We delete $\binom{2n}n*(\beta )\thicksim 3$ and $\binom{2n}n*(\vartheta
)=10^{*}\thicksim A*(e).$

\[
\]
\textbf{3.2 Yifan Yang pullbacks of hypergeometric 5-th order differential
equations.}

Consider the 5-th order hypergeometric differential equation 
\[
L_5=\theta ^5-cx(\theta +a_1)(\theta +a_2)(\theta +a_3)(\theta +a_4)(\theta
+a_5) 
\]

\textbf{Proposition 1.}

If 
\[
a_1=\frac 12 
\]
\[
a_3=1-a_2 
\]
\[
a_5=1-a_4 
\]
then $L_5$ satisfies the C-Y condition

\textbf{Proof: }Write $L_5$ as 
\[
w^{(5)}+b_4w^{(4)}+b_3w^{(3)}+b_2w^{\prime \prime }+b_1w^{\prime }+b_0w=0 
\]
where 
\[
b_4=\frac 52\cdot \frac{4-5cx}{x(1-cx)} 
\]
\[
b_3=\frac{25-c(42+a_2+a_4-a_2^2-a_4^2)x}{x^2(1-cx)} 
\]
\[
b_2=\frac 32\cdot \frac{10-c(26+3a_2+3a_4-3a_2^2-3a_4^2)x}{x^3(1-cx)} 
\]
satisfy C-Y, i.e. 
\[
b_2=\frac 32b_3^{\prime }+\frac 35b_3b_4-b_4^{\prime \prime }-\frac
65b_4b_4^{\prime }-\frac 4{25}b_4^3 
\]
\[
\]

\textbf{Proposition 2.}

$L_5$ has the solutions 
\[
w_0=\sum_{n=0}^\infty \frac{(a_1)_n(a_2)_n(a_3)_n(a_4)_n(a_5)_n}{n!^5}(cx)^n 
\]
\[
w_1=w_0\log (x) 
\]
\[
+\sum_{n=1}^\infty \frac{(cx)^n}{n!^5}\prod_{j=1}^5(a_j)_n\left\{
-5H_n+\sum_{j=1}^5(\psi (n+a_j)-\psi (a_j))\right\} 
\]
where 
\[
H_n=\sum_{j=1}^n\frac 1j 
\]
\[
\psi (a)=\frac{\Gamma ^{\prime }(a)}{\Gamma (a)} 
\]
\[
\]
Observe that 
\[
\psi (n+a)-\psi (a)=\sum_{j=0}^{n-1}\frac 1{a+j} 
\]
\[
(a)_n^{\prime }=(a)_n(\psi (n+a)-\psi (a))=\sum_{j=0}^{n-1}a(a+1)...\widehat{%
(a+j)}...(a+n-1) 
\]
Using Yinfan Yang's definition of the pullback we get the following result

\textbf{Theorem }

The differential equation 
\[
\theta ^4-cx\left\{ 2(\theta +\frac 12)^4+\frac 12(\frac 72-\alpha ^2-\beta
^2)(\theta +\frac 12)^2+\frac 1{16}-\frac 14(\alpha ^2+\frac 14)(\beta
^2+\frac 14)\right\} 
\]
\[
+c^2x^2(\theta +1+\frac{\alpha +\beta }2)(\theta +1-\frac{\alpha +\beta }%
2)(\theta +1+\frac{\alpha -\beta }2)(\theta +1-\frac{\alpha -\beta }2) 
\]
has the solution 
\[
y_0=\frac 1{\sqrt{1-cx}}\left\{ 
\begin{array}{c}
\sum_{n=0}^\infty (cx)^n\sum_{k=0}^n\frac{\prod_{j=1}^5(a_j)_k(a_j)_{n-k}}{%
k!^5(n-k)!^5}\cdot \\ 
\lbrack 1+(2k-n)\{-5H_k+\sum_{j=1}^5(\psi (k+a_j)-\psi (a_j))\}]%
\end{array}
\right\} ^{1/2} 
\]
where 
\[
a_1=\frac 12 
\]
\[
a_2=\frac 12+\alpha 
\]
\[
a_3=\frac 12-\alpha 
\]
\[
a_4=\frac 12+\beta 
\]
\[
a_5=\frac 12-\beta 
\]
\[
\]

We get the following 14 cases giving CY equations 
\[
\]
\[
\begin{tabular}{|c|c|c|c|c|c|}
\hline
Name & $a_2$ & $a_4$ & c & $\alpha $ & $\beta $ \\ \hline
$\widetilde{1}$ & 1/5 & 2/5 & 4$\cdot 5^5$ & -3/10 & -1/10 \\ \hline
$\widetilde{2}$ & 1/10 & 3/10 & 4$\cdot 8\cdot 10^5$ & -2/5 & -1/5 \\ \hline
$\widetilde{3}$ & 1/2 & 1/2 & 4$\cdot 256$ & 0 & 0 \\ \hline
$\widetilde{4}$ & 1/3 & 1/3 & 4$\cdot 3^6$ & -1/6 & -1/6 \\ \hline
$\widetilde{5}$ & 1/2 & 1/3 & 4$\cdot 432$ & 0 & -1/6 \\ \hline
$\widetilde{6}$ & 1/2 & 1/4 & 4$\cdot 2^{10}$ & 0 & -1/4 \\ \hline
$\widetilde{7}$ & 1/8 & 3/8 & 4$\cdot 2^{16}$ & -3/8 & -1/8 \\ \hline
$\widetilde{8}$ & 1/6 & 1/3 & 4$\cdot 11664$ & -1/3 & -1/6 \\ \hline
$\widetilde{9}$ & 1/12 & 5/12 & 4$\cdot 12^6$ & -5/12 & -1/12 \\ \hline
$\widetilde{10}$ & 1/4 & 1/4 & 4$\cdot 2^{12}$ & -1/4 & -1/4 \\ \hline
$\widetilde{11}$ & 1/4 & 1/3 & 4$\cdot 12^3$ & -1/4 & -1/6 \\ \hline
$\widetilde{12}$ & 1/6 & 1/4 & 4$\cdot 2^{10}3^3$ & -1/3 & -1/4 \\ \hline
$\widetilde{13}$ & 1/6 & 1/6 & 4$\cdot 2^83^6$ & -1/3 & -1/3 \\ \hline
$\widetilde{14}$ & 1/2 & 1/6 & 4$\cdot 2^83^3$ & 0 & -1/3 \\ \hline
\end{tabular}
\]
\[
\]

We have the following table of the differential equations

$\widetilde{\mathbf{1}}$%
\[
\theta ^4-10x(2500\theta ^4+5000\theta ^3+5875\theta ^2+3375\theta +738) 
\]
\[
+62500z^2(5\theta +4)(5\theta +6)(10\theta +9)(10\theta +11) 
\]

$\widetilde{\mathbf{2}}$%
\[
\theta ^4-8x(80000\theta ^4+160000\theta ^3+186000\theta ^2+106000\theta
+22811) 
\]
\[
+2^{10}10^6z^2(10\theta +7)(10\theta +9)(10\theta +11)(10\theta +13) 
\]

$\widetilde{\mathbf{3}}$%
\[
\theta ^4-16x(128\theta ^4+256\theta ^3+304\theta ^2+176\theta
+39)+2^{20}x^2(\theta +1)^4 
\]

$\widetilde{\mathbf{4}}$%
\[
\theta ^4-18x(324\theta ^4+648\theta ^3+765\theta ^2+441\theta +97) 
\]
\[
+2^23^{10}x^2(\theta +1)^2(6\theta +5)(6\theta +7) 
\]

$\widetilde{\mathbf{5}}$%
\[
\theta ^4-12x(288\theta ^4+576\theta ^3+682\theta ^2+394\theta +87) 
\]
\[
+144x^2(12\theta +11)^2(12\theta +13)^2 
\]

$\widetilde{\mathbf{6}}$%
\[
\theta ^4-16x(512\theta ^4+1024\theta ^3+1208\theta ^2+696\theta +153) 
\]
\[
+4096x^2(8\theta +7)^2(8\theta +9)^2 
\]

$\widetilde{\mathbf{7}}$%
\[
\theta ^4-16x(32768\theta ^4+65536\theta ^3+76544\theta ^2+43776\theta
+9495) 
\]
\[
+2^{26}z^2(4\theta +3)(4\theta +5)(8\theta +7)(8\theta +9) 
\]

$\widetilde{\mathbf{8}}$%
\[
\theta ^4-36x(2592\theta ^4+5184\theta ^3+6066\theta ^2+3474\theta +755) 
\]
\[
+2^43^{10}z^2(4\theta +3)(4\theta +5)(12\theta +11)(12\theta +13) 
\]

$\widetilde{\mathbf{9}}$%
\[
\theta ^4-144x(165888\theta ^4+331776\theta ^3+386496\theta ^2+220608\theta
+47711) 
\]
\[
+2^{22}3^{10}z^2(4\theta +3)(4\theta +5)(6\theta +5)(6\theta +7) 
\]

$\widetilde{\mathbf{10}}$%
\[
\theta ^4-16x(2048\theta ^4+4096\theta ^3+4800\theta ^2+2752\theta +599) 
\]
\[
+2^{24}z^2(\theta +1)^2(4\theta +3)(4\theta +5) 
\]

$\widetilde{\mathbf{11}}$%
\[
\theta ^4-12x(1152\theta ^4+2304\theta ^3+2710\theta ^2+1558\theta +341) 
\]
\[
+144z^2(24\theta +19)(24\theta +23)(24\theta +25)(24\theta +29) 
\]

$\widetilde{\mathbf{12}}$%
\[
\theta ^4-48x(4608\theta ^4+9216\theta ^3+10744\theta ^2+6136\theta +1325) 
\]
\[
+2^{12}3^2z^2(24\theta +17)(24\theta +23)(24\theta +25)(24\theta +31) 
\]

$\widetilde{\mathbf{13}}$%
\[
\theta ^4-12^2x(10368\theta ^4+20736\theta ^3+24048\theta ^2+13680\theta
+2927) 
\]
\[
+12^{10}z^2(\theta +1)^2(3\theta +2)(3\theta +4) 
\]

$\widetilde{\mathbf{14}}$%
\[
\theta ^4-48x(1152\theta ^4+2304\theta ^3+2704\theta ^2+1552\theta +339) 
\]
\[
+2^{16}3^2z^2(6\theta +5)^2(6\theta +7)^2 
\]
\[
\]

We have the equivalences $\widehat{m}\sim \widetilde{m}$ for m=1,2,...14.
Some of the Hadamard products are also equivalent to certain $\widetilde{m,}$
due to the following table $\theta $%
\[
\begin{tabular}{|c|c|c|c|c|}
\hline
{*} & (e) & (h) & (i) & (j) \\ \hline
(e) & $\widetilde{3}$ & $\widetilde{5}$ & $\widetilde{6}$ & $\widetilde{14}$
\\ \hline
(h) &  & $\widetilde{4}$ & $\widetilde{11}$ & $\widetilde{8}$ \\ \hline
(i) &  &  & $\widetilde{10}$ & $\widetilde{12}$ \\ \hline
(j) &  &  &  & $\widetilde{13}$ \\ \hline
\end{tabular}
\]
\[
\]
Let $y_0(x)$ be the solution of $(e)*(e)$ and $\widetilde{y}_0(x)$ of $%
\widetilde{3}.$ Then we have the transformation 
\[
\widetilde{y}_0(x)=f(x)y_0(g(x)) 
\]
where 
\[
f(x)=1+480x+383488x^2+330493952x^3+... 
\]
\[
g(x)=x+128x^2+91920x^3+52555776x^4+... 
\]
We have not found any explicit formulas for $f$ and $g.$ On the contrary, if 
$\widehat{y}_0(x)$ is the solution of $\widehat{3}$ then 
\[
y_0(x)=\frac 1{\sqrt{1-2^{16}x^2}}\widehat{y}_0(\frac x{(1+2^8x)^2}) 
\]
\[
\]

We also notice the following equivalences 
\[
\begin{tabular}{|c|c|c|c|c|}
\hline
{*} & $\beta $ & $\iota $ & $\theta $ & $\varkappa $ \\ \hline
A & A*(e) & B*(e) & C*(e) & D*(e) \\ \hline
B & A*(h) & B*(h) & C*(h) & D*(h) \\ \hline
C & A*(i) & B*(i) & C*(i) & D*(i) \\ \hline
D & A*(j) & B*(j) & C*(j) & D*(j) \\ \hline
\end{tabular}
\]
\[
\]

Here $\beta ,\iota ,\theta ,\varkappa $ refer to third degree differential
equations in [2] and e.g. the upper most left entry means that the pullback
of the 5-th order equation $A*\beta $ is equivalent to the 4-th order,
degree 2 equation $A*(e)$

Hence we have replaced 14+10+16=40 higher degree equations in the table with
equations of degree 2. 
\[
\]

\textbf{3.3. Sporadic CY equations of degree 2.}

The following CY equations were found in a computer search by van Enckevort
and van Straten (see [3] and[5]. They found also \#181=$\widetilde{4}$ ,
\#204=$\widetilde{3}$ and \#31$\thicksim 3.$%
\[
\]

\textbf{18.} 
\[
\theta ^4-4x(2\theta +1)^2(3\theta ^2+3\theta +1)-16x^2(2\theta +1)(2\theta
+3)(4\theta +3)(4\theta +5) 
\]
\[
\]

\textbf{26.} 
\[
\theta ^4-2x(2\theta +1)^2(13\theta ^2+13\theta +4)-12x^2(2\theta
+1)(2\theta +3)(3\theta +2)(3\theta +4) 
\]
\[
\]

\textbf{28.} 
\[
\theta ^4-x(65\theta ^4+130\theta ^3+105\theta ^2+40\theta +6)+4x^2(\theta
+1)^2(4\theta +3)(4\theta +5) 
\]
\[
\]

\textbf{84.} 
\[
\theta ^4-4x(32\theta ^4+64\theta ^3+63\theta ^2+31\theta +6)+256x^2(\theta
+1)^2(4\theta +3)(4\theta +5) 
\]
\[
\]

\textbf{182.} 
\[
\theta ^4-x(43\theta ^4+86\theta ^3+77\theta ^2+34\theta +6)+12x^2(\theta
+1)^2(6\theta +5)(6\theta +7) 
\]
\[
\]

\textbf{183.} 
\[
\theta ^4-4x(2\theta +1)^2(7\theta ^2+7\theta +3)+48x^2(2\theta +1)(2\theta
+3)(4\theta +3)(4\theta +5) 
\]
\[
\]

\textbf{205.} 
\[
\theta ^4-x(59\theta ^4+118\theta ^3+105\theta ^2+46\theta +8)+96x^2(\theta
+1)^2(3\theta +2)(3\theta +4) 
\]
\[
\]

We deleted \#31 which is equivalent to \#3 which is of degree one.

There are a number of 4-th order differential equations that fulfill all the
conditions of CY equations but they are factorable. Most of them were found
by van Enckevort and van Straten. We list them here hoping that somebody can
give a geometric explanation of their existence. 
\[
\]

\textbf{(i)} 
\[
\theta ^4-6x(36\theta ^4+72\theta ^3+83\theta ^2+47\theta +10)+36x^2(3\theta
+2)(3\theta +4)(6\theta +5)(6\theta +7) 
\]
\[
\]

\textbf{(ii)} 
\[
\theta ^4-4x(32\theta ^4+64\theta ^3+74\theta ^2+42\theta +9)+16x^2(4\theta
+3)^2(4\theta +5^2) 
\]
\[
\]

\textbf{(iii)} 
\[
\theta ^4-4x(128\theta ^4+256\theta ^3+294\theta ^2+166\theta
+35)+16x^2(8\theta +5)(8\theta +7)(8\theta +9)(8\theta +11) 
\]
\[
\]

\textbf{(iv)} 
\[
\theta ^4-12x(6\theta +1)(6\theta +5)(2\theta ^2+2\theta +1)+144x^2(6\theta
+1)(6\theta +5)(6\theta +7)(6\theta +11) 
\]
\[
\]

\textbf{(v)} 
\[
\theta ^4-3x(3\theta +1)(3\theta +2)(2\theta ^2+2\theta +1)+9x^2(3\theta
+1)(3\theta +2)(3\theta +4)(3\theta +5) 
\]
\[
\]

\textbf{(vi)} 
\[
\theta ^4-4x(4\theta +1)(4\theta +3)(2\theta ^2+2\theta +1)+16x^2(4\theta
+1)(4\theta +3)(4\theta +5)(4\theta +7) 
\]
\[
\]

\textbf{(vii)} 
\[
\theta ^4-4x(2\theta +1)^2(2\theta ^2+2\theta +1)+16x^2(2\theta
+1)^2(2\theta +3)^2 
\]
\[
\]

\textbf{(viii)} 
\[
\theta ^4-2x(2\theta +1)^2(\theta ^2+\theta +1)+4x^2(\theta +1)^2(2\theta
+1)(2\theta +3) 
\]
\[
\]

\textbf{3.4. CY equations of degree 3.}

So far we have only found the following CY equations of degree $3.$ The
notation $\widetilde{214}$ means an equation which is equivalent to $214.$%
\[
\]

\textbf{34.}(Verrill) 
\[
\theta ^4-x(35\theta ^4+70\theta ^3+63\theta ^2+28\theta +5) 
\]
\[
+x^2(\theta +1)^2(259\theta ^2+518\theta +285)-225x^3(\theta +1)^2(\theta
+2)^2 
\]
\[
\]
\[
A_n=\sum_{i+j+k+l+m=n}\Big(\frac{n!}{i!j!k!l!m!}\Big)^2 
\]
\[
\]

$\widetilde{\mathbf{145}}.$%
\[
\theta ^4+2\cdot 3^2x(486\theta ^4+324\theta ^3+279\theta ^2+117\theta +16) 
\]
\[
+2^23^8x^2(972\theta ^4+1296\theta ^3+1035\theta ^2+306\theta +32) 
\]
\[
+2^43^{16}x^3(2\theta +1)^2(3\theta +1)83\theta +2) 
\]
\[
\]
\[
A_n=27^n\binom{2n}n^2\binom{3n}n\sum_k(-1)^{n+k}\binom nk\binom{-1/3}k\binom{%
-2/3}k\binom{n+k}n^{-1} 
\]
\[
\]

$\widetilde{\mathbf{155}}.$%
\[
\theta ^4+2^4x(3072\theta ^4+2048\theta ^3+1728\theta ^2+704\theta +87) 
\]
\[
+2^{18}x^2(3072\theta ^4+4096\theta ^3+3200\theta ^2+896\theta +87) 
\]
\[
+2^{36}x^3(2\theta +1)^2(4\theta +1)(4\theta +3) 
\]
\[
\]
\[
A_n=64^n\binom{2n}n^2\binom{4n}{2n}\sum_k(-1)^{n+k}\binom nk\binom{-1/4}k%
\binom{-3/4}k\binom{n+k}n^{-1} 
\]
\[
\]

$\widetilde{\mathbf{165.}}$%
\[
\theta ^4-3^2x(33\theta ^4+66\theta ^3+57\theta ^2+24\theta +4) 
\]
\[
+2^33^6x^2(\theta +1)^2(5\theta ^2+10\theta +4) 
\]
\[
-2^23^{10}x^3(\theta +1)(\theta +2)82\theta +1)82\theta +5) 
\]
\[
\]
\[
A_n=27^n\binom{2n}n\sum_k(-1)^{n+k}\binom nk\binom{3k}n\binom{-1/3}k\binom{%
-2/3}k 
\]
\[
\]
\[
\]

$\widetilde{\mathbf{214}}.$%
\[
\theta ^4-2x(102\theta ^4+204\theta ^3+155\theta ^2+53\theta +7) 
\]
\[
+2^2x^2(\theta +1)^2(396\theta ^2+792\theta +311) 
\]
\[
-2^47^2x^3(\theta +1)8\theta +2)82\theta +1)(2\theta +5) 
\]
\[
\]
\[
A_n=\binom{2n}n\sum_{i,j}(-1)^{i+j}\binom ni\binom nj\binom{i+j}j^3 
\]
\[
\]

\textbf{227.} 
\[
\theta ^4-2^23^2x(132\theta ^4+264\theta ^3+201\theta ^2+69\theta +10) 
\]
\[
+2^93^6x^2(20\theta ^4+80\theta ^3+107\theta ^2+54\theta +10) 
\]
\[
+2^{12}3^{10}x^3(2\theta +1)^2(2\theta +3)^2 
\]
\[
\]
\[
A_n=432^n\binom{2n}n\sum_k(-1)^k\binom nk\binom{3k}n\binom{-1/6}k\binom{-5/6}%
k 
\]
\[
\]

\textbf{228.} 
\[
\theta ^4-2^2x(176\theta ^4+352\theta ^3+289\theta ^2+113\theta +18) 
\]
\[
+2^{11}x^2(80\theta ^4+320\theta ^3+449\theta ^2+258\theta +54) 
\]
\[
-3\cdot 2^{16}x^3(2\theta +1)(2\theta +5)(4\theta +3)(4\theta +9) 
\]
\[
\]
\[
A_n=64^n\binom{2n}n\sum_k(-1)^{n+k}\binom nk\binom{3k}n\binom{-1/4}k\binom{%
-3/4}k 
\]
\[
\]

\textbf{4. Difference equations of low order.}

Most of the CY equations were found by using Maple's ''Zeilberger'' on a
simple sum of binomial coefficients. Maple answers with a recursion formula
of smallest order, but does not care if the degree is high. Since when
converting to differential equations, the degree becomes the order of the
differential equation, we often end up with order larger than 4 (usually 6).
In most cases the differential equation $L$ can be factored $L=L_1L_2$ where
the order of $L_2$ is 4.

Consider the noncommutative \textsl{algebra of difference operators } 
\[
\mathbf{C}\left\langle n,N\right\rangle 
\]
where 
\[
Nf(n)=f(n+1) 
\]
We have the relation 
\[
Nn=nN+N 
\]
We have the \textsl{Weyl algebra }of differential operators with polynomial
coefficients 
\[
\mathbf{C}\left\langle x,\theta \right\rangle 
\]
where $\theta =x\frac d{dx}$ with the relation 
\[
\theta x=x\theta +x 
\]

\textbf{Proposition.}

There is an isomorphism between the differerence operator algebra 
\[
\mathbf{C}\left\langle n,N\right\rangle /(Nn-nN-N) 
\]
and the Weyl algebra 
\[
\mathbf{C}\left\langle x,\theta \right\rangle /(\theta x-x\theta -x) 
\]
given by the map 
\[
n\rightarrow -\theta 
\]
\[
N\rightarrow x 
\]
Hereby 
\[
order\rightarrow \deg ree 
\]
\[
\deg ree\rightarrow order 
\]
\[
\]

In order to solve the differential equations we need to write them in
ordinary fashion. For this the following formulas are useful 
\[
\theta ^n=\sum_{j=1}^nS(n,j)x^j\frac{d^j}{dx^j} 
\]
\[
x^n\frac{d^n}{dx^n}=\sum_{j=1}^ns(n,j)\theta ^n 
\]
where $S(n,j)$ and $s(n,j)$ are \textsl{Stirling numbers }of second and
first kind respectively. This gives a very nice way of defining the Stirling
numbers without referring to combinatorics. 
\[
\]

\textbf{Example }(\#232) Let 
\[
A_n=\binom{2n}n^2\sum_k\binom nk^2\binom{3n}{n+k} 
\]
Maple gives the following difference operator 
\[
(n+2)^4Q_0(n+2)N^2 
\]
\begin{eqnarray*}
&&-(37345n^6+272085n^5+806321n^4+123965n^3 \\
&&+1038758n^2+448272n+77824)N
\end{eqnarray*}
\[
-384(2n+1)^2(3n+1)(3n+2)Q_0(n+3) 
\]
where 
\[
Q_0(n)=77n^2-209n+142 
\]
Converting this to a differential operator we get 
\[
L=\theta ^4Q_0(\theta ) 
\]
\[
-x(37345\theta ^6+48015\theta ^5+6071\theta ^4-11683\theta ^3-2944\theta
^2+78\theta +240) 
\]
\[
-384x^2(2\theta +1)^2(3\theta +1)(3\theta +2)Q_0(\theta +3) 
\]
Factorization gives 
\[
L=L_1L_2 
\]
where 
\[
L_1=5^2Q_0(\theta +4)-5x(201509\theta ^2+1453353\theta +2642894) 
\]
\[
+2^63x^2(51821^2+287078\theta +499706) 
\]
\[
+2^{10}3^25x^3(14707\theta ^2+158741\theta +453284) 
\]
\[
-2^{19}3^5x^4Q_0(\theta +7) 
\]
and 
\[
L_2=5^2\theta ^4-5x(2617\theta ^4+4658\theta ^3+3379\theta ^2+1050\theta
+120) 
\]
\[
-2^63x^2(-673\theta ^4+4871\theta ^3+10282\theta ^2+5410\theta +860) 
\]
\[
+2^{10}3^2x^3(955\theta ^4+4320\theta ^3+3477\theta ^2+1020\theta +100) 
\]
\[
-2^{17}3^3x^4(2\theta +1)^2(3\theta +1)(3\theta +2) 
\]
\[
\]

We also give an example of Hadamard squares. Let 
\[
u_0=\sum_{n=0}^\infty E_nx^n 
\]
be the solution of 
\[
\theta ^2-x(a\theta ^2+a\theta +b)-cx^2(\theta +1)^2 
\]
(one of the 10 cases (a),(b),...(j) in [2]). Then the Hadamard square 
\[
y_0=\sum E_n^2x^n=\sum A_nx^n 
\]
satisfies the recursion of order 3 
\[
U(n)(n+3)^4A_{n+3}-U(n+1)\left\{ U(n)U(n+1)+c(n+2)^4\right\} A_{n+2} 
\]
\[
-cU(n)\left\{ U(n)U(n+1)+c(n+2)^4\right\} A_{n+1}+c^3U(n+1)(n+1)^4A_n=0 
\]
where 
\[
U(n)=a(n+1)^2+a(n+1)+b 
\]
This gives a degree 3 differential equation 
\[
Q_0(\theta )\theta ^4-...+c^3x^3Q_0(\theta +4)(\theta +1)^4 
\]
where 
\[
Q_0(n)=a(n-2)^2+a(n-2)+b 
\]
\[
\]

By duality we believe that difference operators of order 2 and 3 are worth a
table which can be used as a ''superseeker'' to recognize CY differential
equations without solving them. We order them with respect to the absolute
value of the discriminant $D$ of $Q_0(n).$ It is worth mentioning that we
have found numerous formulas for the same coefficient of a CY equation. As
an example we have about 40 different formulas for the coefficient $A_n$ of
\#25. For hundreds of such formulas see [4]. 
\[
\]

\textbf{Table. Superseeker II.} 
\[
\]
\[
\begin{tabular}{|c|c|c|c|c|}
\hline
$\left| D\right| $ & D & Q$_0(n)$ & \# & degree \\ \hline
3 & -3 & $7n^2-23n+19$ & 27 & ? \\ \hline
3 & -3 & $7n^2-19n+13$ & 243 & ? \\ \hline
4 & $-2^2$ & $10n^2-26n+17$ & 237,256 & 2 \\ \hline
5 & $5$ & $5n^2-15n+11$ & 253 & 2 \\ \hline
7 & $-7$ & $14n^2-35n+22$ & 241,33 & 2 \\ \hline
12 & $-2^23$ & $12n^2-30n+19$ & 258 & 2 \\ \hline
12 & $2^23$ & $2n^2-2n-1$ & 56 & 3 \\ \hline
12 & $2^23$ & $2n^2-1+n+11$ & 23 & 3 \\ \hline
15 & $-3\cdot 5$ & $51n^2-147n+106$ & 222 & 2 \\ \hline
15 & $-3\cdot 5$ & $15n^2-45n+34$ & 216 & 3 \\ \hline
15 & $-3\cdot 5$ & $20n^2-55n+38$ & 55 & 2 \\ \hline
15 & $-3\cdot 5$ & $24n^2-57n+34$ & 211 & 2 \\ \hline
16 & $-2^4$ & $20n^2-56n+33$ & 119 & 3 \\ \hline
20 & $2^25$ & $4n^2-18n+19$ & 262 & 3 \\ \hline
27 & $-3^3$ & $27n^2-63n+37$ & 239 & 2 \\ \hline
28 & $2^27$ & $6n^2-26n+27$ & 235 & 4 \\ \hline
32 & $-2^5$ & $24n^2-56n+33$ & 265 & 2 \\ \hline
35 & $-5\cdot 7$ & $21n^2-49n+29$ & 71 & 3 \\ \hline
35 & $-5\cdot 7$ & $21n^2-77n+71$ & 21 & 3 \\ \hline
\end{tabular}
\]
\[
\]
\[
\begin{tabular}{|c|c|c|c|c|}
\hline
$\left| D\right| $ & $D$ & $Q_0(n)$ & \# & degree \\ \hline
37 & $37$ & $41n^2-105n+67$ & 300 & 2 \\ \hline
39 & $-3\cdot 13$ & $20n^2-51n+33$ & 223 & 3 \\ \hline
44 & $-2^211$ & $33n^2-88n+59$ & 278 & 2 \\ \hline
44 & $-2^211$ & $44n^2-110n+69$ & 238,288 & 2 \\ \hline
55 & $-5\cdot 11$ & $77n^2-209n+142$ & 232 & 2 \\ \hline
60 & $-2^23\cdot 5$ & $40n^2-90n+51$ & 277 & 2 \\ \hline
60 & $-2^23\cdot 5$ & $48n^2-126n+83$ & 210 & 2 \\ \hline
105 & $3\cdot 5\cdot 7$ & $12n^2-45n+40$ & 242,259 & 3 \\ \hline
135 & $-3^35$ & $27n^2-99n+92$ & 266 & 3 \\ \hline
140 & $-2^25\cdot 5$ & $52n^2-134n+87$ & 282 & 2 \\ \hline
160 & $2^55$ & $8n^2-8n-3$ & 83 & 3 \\ \hline
176 & $2^411$ & $44n^2-88n+43$ & 254,295 & ? \\ \hline
195 & $-3\cdot 5\cdot 13$ & $85n^2-235n+163$ & 99 & 2 \\ \hline
224 & $-2^57$ & $84n^2-196n+115$ & 289 & 2 \\ \hline
231 & $-3\cdot 7\cdot 11$ & $55n^2-143n+94$ & 117,118 & 3 \\ \hline
231 & $-3\cdot 7\cdot 11$ & $55n^2-187n+160$ & 22,212 & 3 \\ \hline
240 & $-2^43\cdot 5$ & $204n^2-432n+229$ & 225 & 2 \\ \hline
252 & $-2^23^27$ & $56n^2-154n+107$ & 215 & 3 \\ \hline
255 & $-3\cdot 5\cdot 17$ & $87n^2-327n+308$ & 279 & 3 \\ \hline
\end{tabular}
\]
\[
\]
\[
\begin{tabular}{|c|c|c|c|c|}
\hline
$\left| D\right| $ & $D$ & $Q_0(n)$ & \# & degree \\ \hline
288 & $2^53^2$ & $8n^2-40n+41$ & 119 & 3 \\ \hline
320 & $-2^65$ & $48n^2-176n+163$ & 246 & 3 \\ \hline
320 & $-2^65$ & $48n^2-112n+67$ & 247 & 3 \\ \hline
345 & $3\cdot 5\cdot 23$ & $92n^2-299n+242$ & 226 & 3 \\ \hline
385 & $5\cdot 7\cdot 11$ & $44n^2-143n+114$ & 219 & 3 \\ \hline
399 & $-3\cdot 7\cdot 19$ & $285n^2-969n+824$ & 59 & 3 \\ \hline
399 & $-3\cdot 7\cdot 19$ & $102n^2-309n+235$ & 218 & 3 \\ \hline
455 & $-5\cdot 7\cdot 13$ & $156n^2-403n+261$ & 109 & 2 \\ \hline
495 & $-3^25\cdot 11$ & $170n^2-415n+254$ & 192 & 2 \\ \hline
495 & $-3^25\cdot 11$ & $88n^2-231n+153$ & 260 & 3 \\ \hline
640 & $2^75$ & $32n^2-64n+27$ & 261 & 3 \\ \hline
1463 & $-7\cdot 11\cdot 19$ & $171n^2-551n+446$ & 198 & ? \\ \hline
1564 & $-2^217\cdot 23$ & $184n^2-414n+235$ & 264 & 2 \\ \hline
1664 & $-2^713$ & $364n^2-780n+419$ & 294 & 2 \\ \hline
2156 & $2^27^211$ & $110n^2-370n+313$ & 217 & 3 \\ \hline
2176 & $2^717$ & $96n^2-176n+75$ & 276 & 2 \\ \hline
2560 & $-2^95$ & $128n^2-416n+343$ & 275 & 3 \\ \hline
2665 & $5\cdot 13\cdot 41$ & $164n^2-533n+429$ & 274 & 3 \\ \hline
3135 & $-3\cdot 5\cdot 11\cdot 19\cdot $ & $154n^2-407n+274$ & 231 & 3 \\ 
\hline
\end{tabular}
\]
\[
\]
\[
\begin{tabular}{|c|c|c|c|c|}
\hline
$\left| D\right| $ & $D$ & $Q_0(n)$ & \# & degree \\ \hline
3335 & $-5\cdot 23\cdot 29$ & $203n^2-551n+378$ & 224 & 3 \\ \hline
4180 & $2^25\cdot 11\cdot 19$ & $231n^2-484n+249$ & 230 & 3 \\ \hline
5831 & $-7^317$ & $340n^2-1037n+795$ & 234 & 3 \\ \hline
9204 & $2^23\cdot 13\cdot 59$ & $295n^2-944n+763$ & 248 & 3 \\ \hline
14400 & $-2^63^25^2$ & $288n^2-936n+773$ & 249 & 3 \\ \hline
17199 & $-3^37^213$ & $572n^2-1599n+1125$ & 297 & 3 \\ \hline
17415 & $-3^45\cdot 43$ & $324n^2-1053n+869$ & 273 & 3 \\ \hline
39767 & $-7\cdot 13\cdot 19\cdot 23$ & $884n^2-2405n+1647$ & 208 & 3 \\ 
\hline
44591 & $-17\cdot 43\cdot 61$ & $946n^2-2623n+1830$ & 209 & 3 \\ \hline
60480 & $2^63^35\cdot 7$ & $432n^2-1584n+1487$ & 268 & 3 \\ \hline
64496 & $-2^429\cdot 139$ & $1740n^2-3712n+1989$ & 305 & 2 \\ \hline
64844 & $-2^213\cdot 29\cdot 43$ & $812n^2-2378n+1761$ & 240 & 2 \\ \hline
104719 & $-23\cdot 29\cdot 157$ & $1595n^2-5191n+4240$ & 19 & 3 \\ \hline
170375 & $5^329\cdot 47$ & $1457n^2-4277n+3168$ & 195 & ? \\ \hline
702075 & $-3\cdot 5^2\cdot 11\cdot 23\cdot 37$ & $1771n^2-4807n+3361$ & 252
& 3 \\ \hline
959040 & $-2^63^45\cdot 37$ & $2592n^2-8424n+6937$ & 272 & 3 \\ \hline
5274751 & $23\cdot 79\cdot 2903$ & $6557n^2-18565n+13342$ & 250 & 3 \\ \hline
\end{tabular}
\]
\[
\]

\textbf{Acknowledgements.}

Firstly I want to thank Christian van Enckevort and Duco van Straten who
found many of the equations listed here. Further I thank Wadim Zudilin who
has inspired much of this work and communicated the pullback of Yinfan to
me. Finally I thank Yinfan Yang , whose pullback caused the writing of this
paper.

\[
\]

\textbf{References.}:

\textbf{1. }G.Almkvist, Str\"{a}ngar i m\aa nsken I, Normat 51 (2003),
22-33, II, Normat 51 (2003), 63-79

\textbf{2. }G.Almkvist, W.Zudilin, Differential equations, mirror maps and
zeta values. In Mirror symmetry V, Proceedings of BIRS workshop on
Calabi-Yau Varieties and Mirror Symmetry, December 6-11, 2003

\textbf{3.} G.Almkvist,C. van Enckevort, D. van Straten, W.Zudilin, Tables
of Calabi-Yau equations, math.AG/0507430

\textbf{4. }G.Almkvist, Some binomial identities related to Calabi-Yau
differential equations, In preparation

\textbf{5. }C.van Enckevort, D.van Strayen, Monodromy calculations for
fourth order equations of Calabi-Yau type. In Mirror Symmetry V, Proceedings
of BIRS workshop on Calabi-Yau Varieties and Mirror Symmetry, December 6-11,
2003.

\textbf{6. }Yinfan Yang, Pullbacks, personal communication

\textbf{7. }D.Zagier, Integral solutions of Ap\'{e}ry-like recurrence
equations, Manuscript 2003 
\[
\]

Math Dept

Univ of Lund

Box 118

22100 Lund, Sweden

gert@maths.lth.se

\end{document}